\documentclass{notices}

\usepackage{amsmath, amssymb}
\usepackage{graphicx,epstopdf,color}

\newcommand{\R}{{\mathbb{R}}}
\newcommand{\N}{{\mathbb{N}}}
\newcommand{\e}{\varepsilon}
\newcommand{\eps}{\varepsilon}
\renewcommand{\epsilon}{\varepsilon}

\renewcommand{\le}{\leqslant}

\renewcommand{\ge}{\geqslant}

\title{A comparison between the nonlocal and the classical worlds:\\
minimal surfaces, phase transitions, and geometric flows}

\author{
  Serena Dipierro
  \affil{
    University of Western Australia, Department of
Mathematics and Statistics, 35 Stirling Highway, WA 6009, Crawley,
Australia. Email address: serena.dipierro@uwa.edu.au
    }
}

\begin{document}

\maketitle

The nonlocal world
presents an abundance of surprises and wonders to discover.
These special properties of the nonlocal world are usually the consequence
of long-range interactions, which, especially in presence of
geometric structures and nonlinear phenomena, end up producing
a variety of novel patterns.
We will briefly discuss some of these features,
focusing on the case of (non)local minimal surfaces, (non)local
phase coexistence models, and (non)local geometric flows.

\section{Nonlocal minimal surfaces}

\begin{figure}
    \centering
    \includegraphics[height=4.5cm]{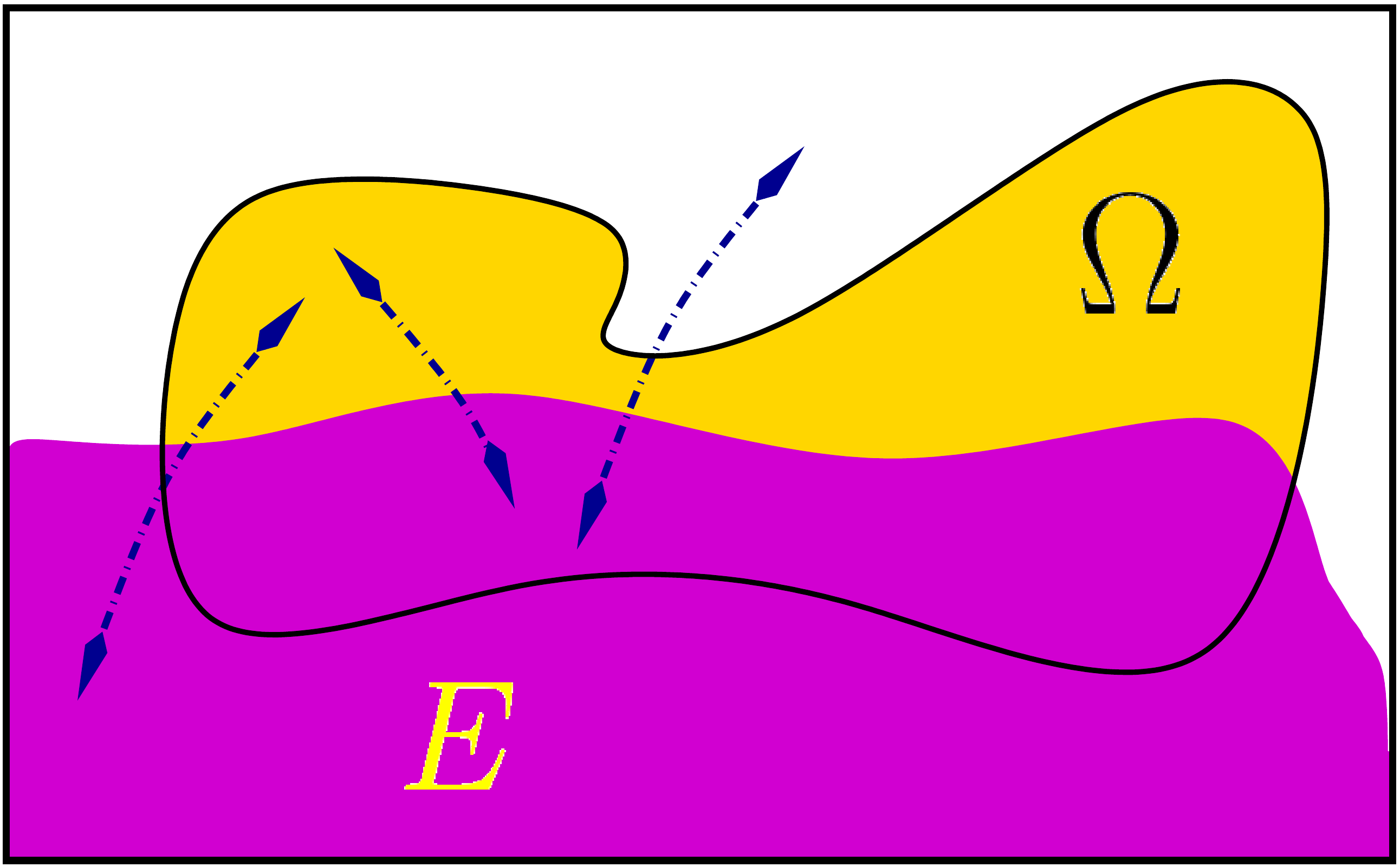}
    \caption{\footnotesize\it Long-range interactions leading to the fractional perimeter of the set~$E$ in~$\Omega$.}
    \label{NESS}
\end{figure}

In~\cite{MR2675483}, a new notion of nonlocal perimeter
has been introduced, and the study of the corresponding minimizers
has started
(related energy functionals had previously arisen in~\cite{MR1111612}
in the context of phase systems and fractals).
The simple, but deep idea, grounding this new definition consists
in considering pointwise interactions between disjoint sets,
modulated by a kernel. The prototype of these
interactions considers kernels which have translational,
rotational, and dilation invariance. Concretely,
given~$\sigma\in(0,1)$, one considers the $\sigma$-interaction
between two disjoint sets~$F$ and~$G$ in~$\R^n$ as defined by
\begin{equation}\label{IN3r45T}
{\mathcal{I}}_\sigma(F,G):=\iint_{F\times G}\frac{dx\,dy}{|x-y|^{n+\sigma}},\end{equation}
and the $\sigma$-perimeter of a set~$E$ in~$\R^n$
as the $\sigma$-interaction between~$E$ and its complement~$E^c$, namely
\begin{equation}\label{PERs} {\rm Per}_\sigma(E,\R^n):=
{\mathcal{I}}_\sigma(E,E^c).\end{equation}
To deal with local minimizers, given a domain~$\Omega\subset\R^n$
(say, with sufficiently smooth boundary), it is also convenient
to introduce the notion of $\sigma$-perimeter of a set~$E$ in~$\Omega$.
To this end, one can consider the interaction in~\eqref{PERs}
as composed by four different terms (by considering
the intersections of~$E$ and~$E^c$ with~$\Omega$ and~$\Omega^c$),
namely we can rewrite~\eqref{PERs} as
\begin{equation}\label{PERs2} \begin{split}&{\rm Per}_\sigma(E,\R^n)=
{\mathcal{I}}_\sigma(E\cap\Omega,E^c\cap\Omega)\\&\quad+
{\mathcal{I}}_\sigma(E\cap\Omega,E^c\cap\Omega^c)+
{\mathcal{I}}_\sigma(E\cap\Omega^c,E^c\cap\Omega)\\&\quad+
{\mathcal{I}}_\sigma(E\cap\Omega^c,E^c\cap\Omega^c).\end{split}\end{equation}
Among the four terms in the right hand side of~\eqref{PERs2},
the first three terms take into account interactions
in which at least one contribution comes from~$\Omega$
(specifically, all the contributions to~${\mathcal{I}}_\sigma(E\cap\Omega,E^c\cap\Omega)$
come from~$\Omega$, while the contributions to~$
{\mathcal{I}}_\sigma(E\cap\Omega,E^c\cap\Omega^c)$ and~$
{\mathcal{I}}_\sigma(E\cap\Omega^c,E^c\cap\Omega)$ come
from the interactions of points in~$\Omega$ with points in~$\Omega^c$).
The last term in~\eqref{PERs2} is structurally different,
since it only takes into account contributions
coming from outside~$\Omega$. It is therefore natural to
define the~$\sigma$-perimeter of a set~$E$ in~$\Omega$
as the collection of all the contributions in~\eqref{PERs2}
that take into account points in~$\Omega$, thus defining
\begin{equation}\label{PERs3} \begin{split}&{\rm Per}_\sigma(E,\Omega):=
{\mathcal{I}}_\sigma(E\cap\Omega,E^c\cap\Omega)\\&\quad+
{\mathcal{I}}_\sigma(E\cap\Omega,E^c\cap\Omega^c)+
{\mathcal{I}}_\sigma(E\cap\Omega^c,E^c\cap\Omega),\end{split}\end{equation}
see Figure~\ref{NESS}

This notion of perimeter recovers the classical
one as~$\sigma\nearrow1$ in various senses, see e.g.~\cite{zbMATH02134074, MR1942130, MR2765717, MR2782803, MR3827804},
and indeed the analogy between the minimizers of~\eqref{PERs3}
with respect to prescribed sets outside~$\Omega$
(named ``nonlocal minimal surfaces''
in~\cite{MR2675483})
and the minimizers of the classical perimeter (the ``classical minimal
surfaces'')
has been widely exploited
in~\cite{MR2675483} as well in a series of subsequent articles,
such as~\cite{MR3107529}. Notwithstanding the
important similarities between the nonlocal and the classical settings,
many striking differences arise, as we will also describe in this note.\medskip

Minimizers, and more generally critical points, of the $\sigma$-perimeter
satisfy (with a suitable, weak or viscosity, meaning,
and in a principal value sense) the equation
\begin{equation}\label{HS} \begin{split}&{\mathcal{H}}^\sigma_E(x):=
\int_{\R^n}\frac{\Xi_E(y)}{|x-y|^{n+\sigma}}\,dy=0,\\&
\qquad{\mbox{ for }}x\in(\partial E)\cap\Omega,\end{split}\end{equation}
being
$$ \Xi_E(y):=\begin{cases}
-1 & {\mbox{ if }}x\in E,\\
1 & {\mbox{ if }}x\in E^c,
\end{cases}$$
see~\cite{MR2675483, LUMA}.
The quantity~${\mathcal{H}}^\sigma_E$
can be seen as a ``nonlocal mean curvature'',
and indeed it approaches in various senses the classical mean
curvature~${\mathcal{H}}^1_E$
(though important differences arise, see~\cite{MR3230079}).
Interestingly, the nonlocal mean curvature
measures the size of~$E$ with respect to its complement,
in an integral fashion weighted by the interaction kernel, in a precise
form given by~\eqref{HS}.
This interpretation is indeed closely related with the classical
mean curvature, in which
the size of~$E$ is measured only ``in the vicinity of the boundary'',
since (up to normalizing constants)
$$ {\mathcal{H}}^1_E(x)=\lim_{\rho\searrow0}\frac{1}{\rho^{n+1}}\,
\int_{B_\rho(x)} \Xi_E(y)\,dy,\qquad{\mbox{ for }}x\in\partial E.$$
The settings in~\eqref{PERs} and~\eqref{HS}
also reveal the ``fractional notion'' of these nonlocal minimal surfaces.
As a matter of fact, if one considers the fractional Gagliardo seminorm
defined, for every~$s\in(0,1)$, by
$$ [f]_{H^s(\R^n)}:=\sqrt{
\iint_{\R^n\times\R^n}\frac{|f(x)-f(y)|^2}{|x-y|^{n+2s}}\,dx\,dy},$$
it follows that
$$ [\Xi_E]_{H^{\sigma/2}(\R^n)}^2=8\,{\rm Per}_\sigma(E).
$$
Moreover, defining the fractional Laplacian as
$$ (-\Delta)^s f(x):=\int_{\R^n}\frac{f(x)-f(y)}{|x-y|^{n+2s}}\,dy,$$
and, for a smooth set~$E$, identifying a.e. the function~$\Xi_E$
with
$$ \widetilde\Xi_E(y):=\begin{cases}
-1 & {\mbox{ if $x$ is in the interior of $E$}},\\
1 & {\mbox{ if $x$ is in the interior of $E^c$}},\\
0& {\mbox{ if $x\in\partial E$}},\end{cases}$$
one sees that, for~$x\in\partial E$, \begin{equation}\begin{split}\label{FA}
{\mathcal{H}}^\sigma_E(x)\,&=\,
\int_{\R^n}\frac{\widetilde\Xi_E(y)}{|x-y|^{n+\sigma}}\,dy
\\&=\,\int_{\R^n}\frac{\widetilde\Xi_E(y)-\widetilde\Xi_E(x)}{|x-y|^{n+\sigma}}\,dy
\\&=\,
-(-\Delta)^{\sigma/2}\widetilde\Xi_E(x).
\end{split}\end{equation}
In this sense, one can relate ``geometric'' equations,
such as~\eqref{HS}, with ``linear'' equations
driven by the fractional Laplacian (as in~\eqref{FA}),
in which however the nonlinear feature of the problem is encoded
by the fact that the equation takes place on the boundary of a set
(once again, however, sharp differences arise between
nonlocal minimal surfaces and solutions of linear equations,
as we will discuss in the rest of this note).
\medskip

{F}rom the considerations above, we see how
the study of nonlocal minimal surfaces is therefore
related to the one of hypersurfaces with vanishing nonlocal mean
curvature, and a direction of research
lies in finding correspondences
between critical points of
nonlocal and classical perimeters. With respect
to this point, we recall that double helicoids possess both
vanishing classical mean
curvature and vanishing nonlocal mean
curvature, as noticed in~\cite{MR3532174}.
Also, nonlocal catenoids
have been constructed in~\cite{MR3798717}
by bifurcation methods from the classical case. Differently
from the local situation, the nonlocal catenoids exhibit \label{CATE}
linear growth at infinity, rather than a logarithmic one,
see also~\cite{cozzi-lincei}.\medskip

The cases of hypersurfaces with prescribed nonlocal mean
curvature and of minimizers of nonlocal perimeter functionals
under periodic conditions have also been considered in the literature,
see e.g.~\cite{MR3485130, MR3744919, MR3770173, MR3881478}.
For completeness, we also mention that
we are restricting here to the notion of nonlocal minimal surfaces
of fractional type, as introduced in~\cite{MR2675483}.
Other very interesting, and structurally different,
notions of nonlocal minimal surfaces have also been
considered in the literature, see~\cite{MR3930619, MR3981295, MR3996039}
and the references therein.

\section{Interior regularity for
minimal surfaces}

A first step to understand the geometric properties of
(both classical and nonlocal) minimal surfaces is
to detect their regularity properties and the possible singularities.
While the regularity theory of classical minimal surfaces
is a rather well established topic, many basic regularity problems
in the nonlocal minimal surfaces setting are open, and they
require brand new ideas to be attacked.

\subsection{Interior regularity for classical
minimal surfaces}\label{INTCLA}

Classical minimal surfaces are smooth up to
dimension~$7$.
This has been established in~\cite{MR178385}
when~$n = 3$, in~\cite{MR200816} when~$ n = 4$
and in~\cite{MR233295} when~$n\le7$.

Also, it was conjectured in~\cite{MR233295} that
this regularity result was optimal in dimension~$7$, suggesting
as a possible counterexample in dimension~$8$ the cone
$$ \{(x, y) \in \R^4\times\R^4 {\mbox{ s.t. }} |x|<|y|\}.$$
This conjecture was indeed positively assessed in~\cite{MR250205},
thus showing that classical minimal surfaces
can develop singularities in dimension~$8$ and higher.\medskip

Besides the smoothness,
the analyticity of minimal surfaces in dimension up to~$7$
has been established
in~\cite{MR0093649, MR170091, MR0179651}. See also~\cite{MR301343, MR171198} for more general results.\medskip

Comprehensive discussions about the regularity
of classical minimal surfaces can be found e.g. in~\cite{MR756417, MR1361175, MR2760441, MR3468252, MR3588123}.

\subsection{Interior regularity for nonlocal
minimal surfaces}

Differently from the classical case, the regularity
theory for nonlocal minimal surfaces constitutes a challenging problem
which is still mostly open.
Till now, a complete result holds only
dimension~$2$, since the smoothness of nonlocal
minimal surfaces has been established in~\cite{MR3090533}.

In higher dimensions, the results in~\cite{MR3107529}
give that nonlocal minimal surfaces are smooth
up to dimension~$7$ as long as~$\sigma$ is sufficiently close to~$1$:
namely, for every~$n\in\N\cap[1,7]$ there exists~$\sigma(n)\in[0,1)$
such that all $\sigma$-nonlocal minimal surfaces in dimension~$n$
are smooth provided that~$\sigma\in(\sigma(n),1)$.
The optimal value of~$\sigma(n)$ is unknown (except when~$n=1,2$,
in which case~$\sigma(n)=0$).

This result strongly relies on the fact that the nonlocal perimeter
approaches the classical perimeter when~$\sigma\nearrow1$,
and therefore one can expect that when~$\sigma$ is close to~$1$
nonlocal minimal surfaces inherit the regularity of
classical minimal surfaces.
\medskip

On the other hand, we remark that it is not possible to obtain
information on the regularity of nonlocal minimal surfaces from the
asymptotics of the fractional perimeter as~$\sigma\searrow0$. 
Indeed, it has been proved in~\cite{MR3007726} that the fractional
perimeter converges, roughly speaking, to the Lebesgue
measure of the set when~$\sigma\searrow0$,
and so in this case minimizers can be as wild as they wish.
\medskip

As a matter of fact, the smoothness obtained in~\cite{MR3090533}
and~\cite{MR3107529} is of~$C^{1,\alpha}$-type, the improvement to~$C^\infty$
has been obtained in~\cite{MR3331523}. 

Indeed, nonlocal minimal surfaces
enjoy an ``improvement of regularity'' from locally Lipschitz
to~$C^\infty$ (more precisely, from locally Lipschitz
to~$C^{1,\alpha}$ for any~$\alpha<\sigma$,
thanks to~\cite{MR3680376}, and from~$C^{1,\alpha}$ for
some~$\alpha>\sigma/2$ to~$C^\infty$,
thanks to~\cite{MR3331523}).
\medskip

It is an open problem
to establish whether smooth
nonlocal minimal surfaces are actually analytic.
It is also open to determine
whether or not singular nonlocal minimal surfaces exist
(a preliminary analysis performed in~\cite{MR3798717}
for symmetric cones may lead to the conjecture that
nonlocal minimal surfaces are smooth up to dimension~$6$, but completely
new
phenomena may arise in dimension~$7$).\medskip

Quantitative versions of the regularity results for nonlocal
minimal surfaces have been
obtained in~\cite{MR3981295}.\medskip

As a matter of fact, in~\cite{MR3981295}
also more general interaction kernels have been considered,
and regularity results have been
obtained for more general critical points than minimizers.
In particular, one can consider stable surfaces
with vanishing nonlocal mean curvature as
critical points with ``nonnegative second variations''
of the functional. More precisely, one
says that~$E$ is ``stable'' for the nonlocal perimeter in~$\Omega$
if~${\rm Per}_\sigma(E,\Omega)<+\infty$
and for every vector field~$ X = X(x, t) \in C^{2}(\R^n\times(-1, 1); \R^n)$,
which is compactly supported in~$\Omega$
and whose integral flow is denoted by~$\Phi^t_X$,
and for every~$\e>0$, there exists~$t_0>0$ such that
\begin{equation*}
\begin{split}
&{\rm Per}_\sigma \big(\Phi^t_X(E)\cup E\big) -
{\rm Per}_\sigma (E)+ \e t^2\ge0\\{\mbox{and }}\quad
&{\rm Per}_\sigma \big(\Phi^t_X(E)\cap E\big) -
{\rm Per}_\sigma (E)+ \e t^2\ge0.\end{split}\end{equation*}
for all~$t\in(-t_0, t_0)$.

Interestingly, the notion of stability in the nonlocal regime provides
stronger information with respect to the classical counterpart:
for instance, if~$E$ is stable for the nonlocal perimeter in~$B_{2R}$,
then
\begin{equation}\label{67:293e}
{\rm Per} (E, B_R)\le C R^{n-1},\end{equation}
where~$C>0$ is a constant depending only on~$n$ and~$\sigma$.

We stress that the left hand side in~\eqref{67:293e}
involves the classical perimeter (not the nonlocal perimeter),
hence an estimate of this type is quite informative also
for minimizers (not only for stable sets). In a sense,
up to now, the perimeter estimate in~\eqref{67:293e}
is the only regularity results known for nonlocal minimal
surfaces (and, more generally, for nonlocal stable surfaces) in any dimension
(differently from~\cite{MR3090533}
and~\cite{MR3107529}, this estimate does not imply the smoothness
of the surface, but only a bound on the perimeter).\medskip

Moreover, the right hand side of~\eqref{67:293e}
is uniform with respect to the external data of~$E$. In particular,
wild data are shown to have a possible impact on the perimeter
of~$E$ near the boundary, but not in the interior
(for instance, nonlocal minimal surfaces (and, more general,
nonlocal
stable surfaces) in~$B_2$ have
a uniformly bounded perimeter in~$B_1$).
This is a remarkable property, heavily relying on the nonlocal
structure of the problem, which has no counterpart in the classical
case. As an example, one can consider a family of parallel hyperplanes,
which is a local minimizer for the perimeter:
since each hyperplane produces a certain perimeter contribution
in~$B_1$, no uniform bound can be attained in this situation. In this sense,
a bound as in~\eqref{67:293e} prevents arbitrary
families of possibly perturbed hyperplanes to be stable surfaces
for the nonlocal perimeter.
\medskip

Concerning the regularity of stable sets for the nonlocal perimeter,
we also mention that half-spaces are the only stable cones in~$ \R^3$, if the fractional parameter~$\sigma$
is sufficiently close to~$1$, as proved in~\cite{2017arXiv171008722C}.

\section{The Dirichlet problem for
minimal surfaces}

The counterpart of the theory of minimal
surfaces consists in finding solutions with graphical structure
for given boundary or exterior data.
In the classical setting, this corresponds to studying
graphs with vanishing mean curvature inside a given domain
with a prescribed Dirichlet datum along the boundary of the
domain. Its nolocal counterpart consists in
studying
graphs with vanishing nonlocal mean curvature inside a given domain
with a prescribed datum outside this
domain. We discuss now similarities and differences
between these two problems.

\subsection{The Dirichlet problem for classical minimal surfaces}

A classical problem in geometric analysis
is to seek hypersurfaces with vanishing
mean curvature and prescribed boundary data. Namely,
given a smooth
domain~$\Omega\subset\R^n$
and a boundary datum~$\varphi\in C(\partial\Omega)$, the problem
is to
find~$u\in C^2(\Omega)\cap C(\overline\Omega)$
that solves the Dirichlet problem
\begin{equation}\label{DPB1}
\begin{cases}
{\rm div}\left(\frac{\nabla u}{\sqrt{1+|\nabla u|^2}}\right)=0 & {\mbox{ in }}\Omega,\\
u=\varphi& {\mbox{ on }}\partial\Omega.\end{cases}\end{equation}
A classical approach to this problem
(often referred to with the name of ``Hilbert-Haar existence theory''
see Chapter~3.1 of~\cite{MR795963}
and the references therein) consists in fixing~$ M > 0$
and using the Ascoli Theorem to minimize the area
functional among functions in
\begin{eqnarray*} X_M&:=&
\big\{u \in C^{0,1}(\overline\Omega) {\mbox{ with }}
u=\varphi {\mbox{ on }}\partial\Omega \\&&\qquad\qquad{\mbox{ and }}
[u]_{C^{0,1}(\overline\Omega)}\le M
\big\},\end{eqnarray*}
where, as customary, we denote the Lipschitz seminorm of~$u$
by
$$ [u]_{C^{0,1}(\overline\Omega)}:=\sup_{{x,y\in\Omega
}\atop{x\ne y}}\frac{|u(x)-u(y)|}{|x-y|}.$$
To use this direct minimization approach, one needs
of course to check that~$X_M\ne\varnothing$.
Furthermore, in order to obtain a solution of~\eqref{DPB1},
it is crucial that the Lipschitz seminorm of the minimizer in~$X_M$
is in fact {\em strictly smaller} than~$M$ (as long as~$M$
is chosen conveniently large), so to obtain an
{\em interior} minimum of the area functional, and thus
find~\eqref{DPB1} as the Euler-Lagrange
equation of this minimization procedure
(the Lipschitz bound permitting the use
of uniformly elliptic regularity theory for PDEs,
leading to the desired smoothness of the solution inside the domain).\medskip

Finding sufficient and necessary conditions
for this procedure to work, and, in general,
for obtaining solutions of~\eqref{DPB1} has been
a classical topic of investigation. 
The main lines of this research took into
account a ``bounded slope condition'' on the domain and
the datum that allows one to exploit affine functions as barriers.
In this, the convexity of~$\Omega$
played an important role for the explicit construction of linear
barriers. As a matter of fact,
the first existence results for problem~\eqref{DPB1}
dealt with the planar case~$n=2$ and
convex domains~$ \Omega$, see~\cite{MR1509123, MR1512358, MR1545197}.

Existence results in higher dimensions for convex domains
have been established in~\cite{MR0146506, MR155209}.

The optimal conditions for existence results were discovered
in~\cite{MR222467} and rely on the notion of ``mean convexity''
of the domain: namely, 
if~$\Omega$ has~$C^2$ boundary,
problem~\eqref{DPB1} is solvable for every
continuous boundary datum~$\varphi$ if and only if
the mean curvature of~$\Omega$ is nonnegative
(when~$n=2$, the notion of mean convexity
boils down to the usual convexity).
For general results of this flavor, see
Theorem~16.11 in~\cite{MR1814364}
and the references therein.

\subsection{The Dirichlet problem for nonlocal minimal surfaces}

Given a smooth domain~$\Omega\subset\R^n$,
one considers the cylinder~$\Omega^\star:=\Omega\times\R$
and looks for local minimizers of the nonlocal perimeter
among all the sets with prescribed datum outside~$\Omega^\star$
(more precisely, one seeks
local minimizers of the nonlocal perimeter in~$\widetilde\Omega$,
for every smooth and bounded~$\widetilde\Omega\subset\Omega^\star$,
see~\cite{MR3827804} for all the details of this construction).\medskip

In~\cite{MR3516886}, it is shown that this problem
is solvable in the class of graphs. More precisely,
if~$E$ is a local minimizer of the $\sigma$-perimeter
in~$\Omega^\star$ such that its datum outside~$\Omega^\star$
has a graphical structure, namely
\begin{equation} \label{ED1}
E\setminus\Omega^\star = \{x_{n+1} < u_0(x),\;\, x\in\R^n\setminus\Omega\},\end{equation}
for some continuous function~$ u_0: \R^n \to\R$,
then~$E$ has a graphical structure inside~$\Omega^\star$ as well,
that is
\begin{equation} \label{ED2}
E\cap\Omega^\star = \{x_{n+1} < u(x),
\;\, x\in\Omega\}.\end{equation}
for some continuous function~$ u : \R^n\to\R$.\medskip

\begin{figure}
    \centering
    \includegraphics[height=4.5cm]{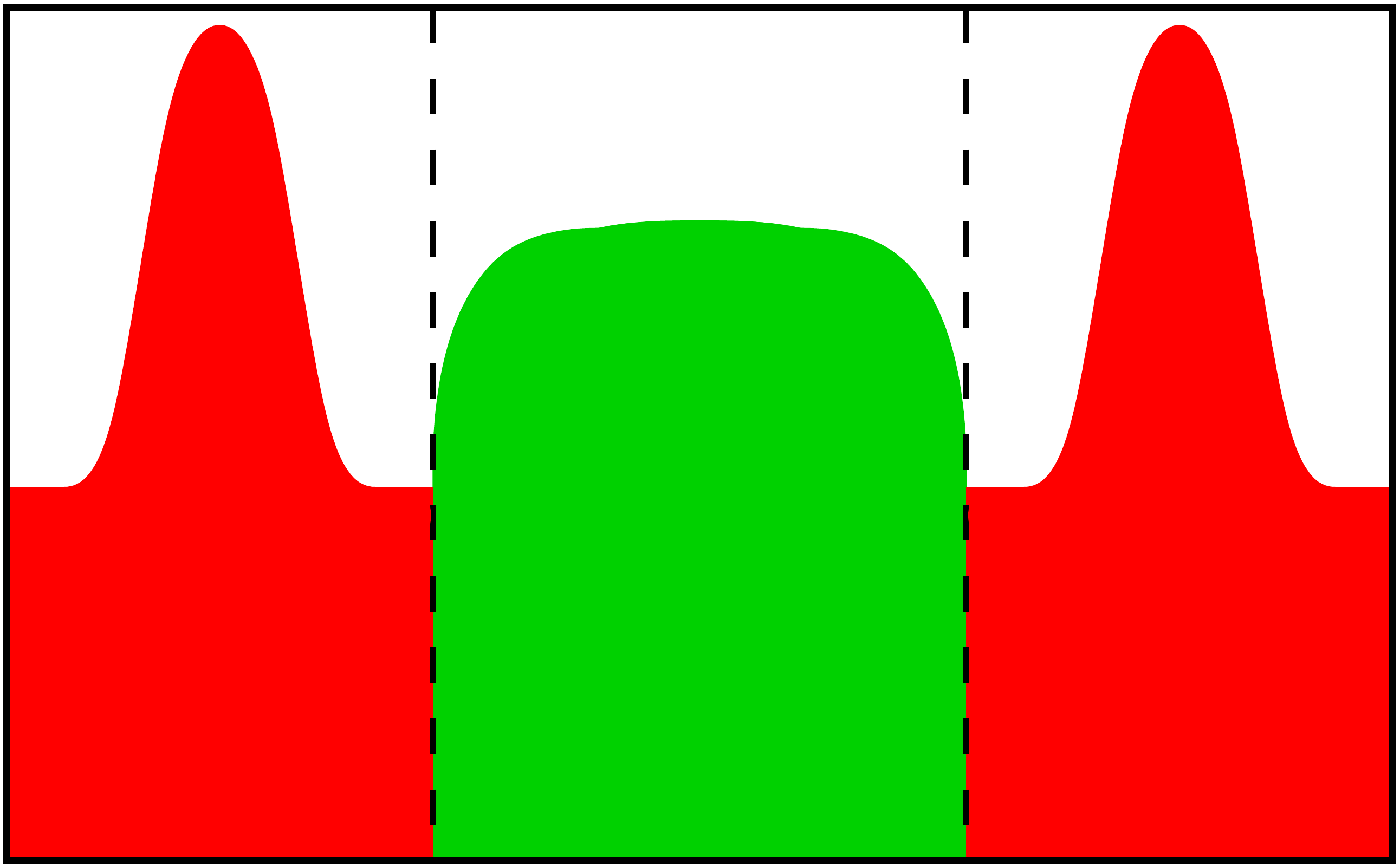}
    \caption{\footnotesize\it An example of stickiness (the red is the prescribed set,
the green is the $s$-minimizer).}
    \label{FF1}
\end{figure}

\begin{figure}
    \centering
    \includegraphics[height=4.5cm]{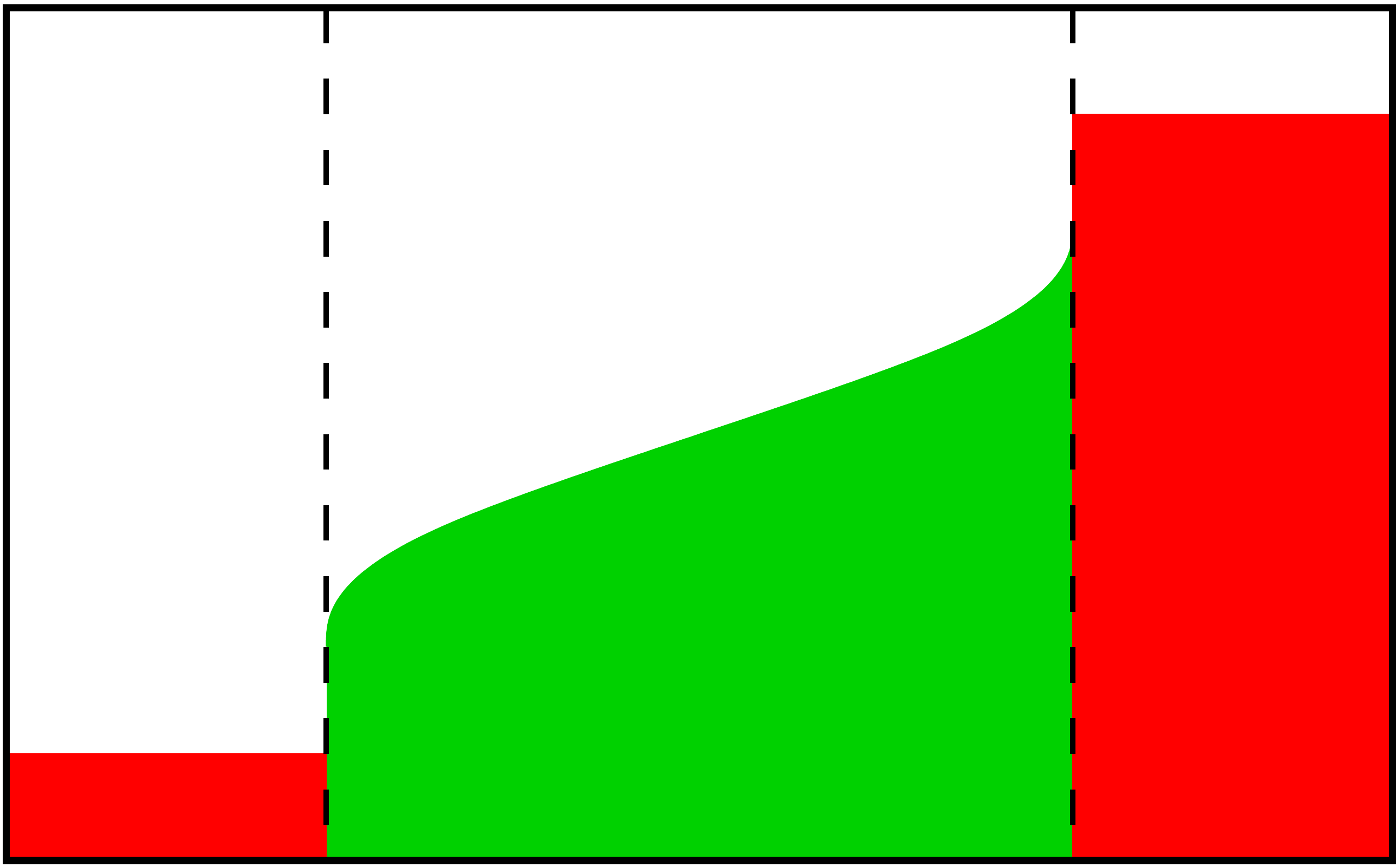}
    \caption{\footnotesize\it Another example of stickiness (the red is the prescribed set,
the green is the $s$-minimizer).}
    \label{FF2}
\end{figure}

Remarkably, in general this problem may lose continuity
at the boundary of~$\Omega^\star$. That is, in the setting
of~\eqref{ED1} and~\eqref{ED2}, it may happen that
\begin{equation}\label{ED3}
\lim_{{x\to\partial\Omega}\atop{x\in\Omega}}u(x)\ne
\lim_{{x\to\partial\Omega}\atop{x\in\Omega^c}}u_0(x).
\end{equation}
The first example of this quite surprising phenomenon was given
in~\cite{MR3596708}, and this is indeed part of a general
and remarkable structure of nonlocal minimal surfaces that
we named ``stickiness'': namely, nonlocal minimal
surfaces have the tendency to stick at the boundary of
the domain (even when the domain is convex,
in sharp contrast with the pattern exhibited by classical minimal surfaces).
See e.g. Figures~\ref{FF1} and~\ref{FF2}
for some qualitative examples.\medskip

Recently, in~\cite{2019arXiv190405393D},
it was discovered that the boundary discontinuity pointed
out in~\eqref{ED3} is indeed a ``generic'' phenomenon
at least in the plane: namely if a given graphical external
datum produces a continuous minimizer, then
an arbitrarily small perturbation of it will produce
a minimizer exhibiting the boundary discontinuity.\medskip

Interestingly, one of the main ingredients in the proof
of the genericity of the boundary discontinuity in~\cite{2019arXiv190405393D}
consists in an ``enhanced boundary regularity'':
namely, nonlocal minimal surfaces in the plane
with a graphical structure that are continuous
at the boundary are automatically~$C^{1,\frac{1+\sigma}2}$
at the  boundary (that is, boundary continuity implies
boundary differentiability).\medskip

A byproduct of this construction, taking into account also the
boundary properties obtained in~\cite{MR3532394},
is also that nonlocal minimal surfaces with a graphical structure
in the plane exhibit a ``butterfly effect'' for their boundary derivatives:
namely, while a nonlocal minimal graph that is continuous
at the boundary has also a finite derivative there,
an arbitrarily small perturbation of the exterior graph will
produce not only a jump at the boundary but also an infinite boundary
derivative (that is, the boundary derivative
switches from a finite, possibly zero, value to infinity only due
to an arbitrarily small perturbation of the datum).
We refer to~\cite{2019arXiv191205794D}
for a survey about this type of
phenomena.\medskip

It is worth recalling that the boundary discontinuity
of nonlocal minimal surfaces is not only a special feature
of the nonlocal world, when compared to the case of classical
minimal surfaces, but also a special feature of the nonlinear
equation of nonlocal mean curvature type, when compared to the case
of linear fractional equations: indeed, as detailed
in~\cite{MR3168912}, solutions of linear equations of the type~$(-\Delta)^s u=f$
in a domain~$\Omega$ with an external datum are typically
H\"older continuous up to the boundary of~$\Omega$,
and no better than this (instead, nonlocal minimal surfaces
are typically discontinuous at the boundary, but when they are
continuous they are also differentiable).
\medskip

Though this discrepancy between nonlocal minimal surfaces
and solutions of linear equations seems a bit surprising at a first
glance, a possible explanation for it lies in the ``accuracy''
of the approximation that a linear equation can offer to an equation
of geometric type, such as the one related to the nonlocal mean curvature.
Namely, linear equations end up being a good approximation
of the geometric equation
with respect to the normal of the surface itself: hence, when
the surface is ``almost horizontal'' the graphical structure of
the geometric equation is well shadowed by its linear counterpart, but
when the surface is ``almost vertical'' the linear counterpart
of the geometric equation should take into account not the graph of
the original surface but its ``inverse'' (i.e., the graph describing
the surface in the horizontal, rather than vertical, direction).
\medskip

We also refer to~\cite{2019arXiv190701498D} for
a first analysis of the stickiness phenomenon in dimension~$3$.
\medskip

Many examples of nonlocal minimal surfaces
exhibiting the stickiness phenomenon are studied in detail
in~\cite{MR3926519}, also with respect to a suitable parameter
measuring the weighted 
measure of the datum at infinity.
\medskip

In addition, very accurate and fascinating
numerical simulations on the stickiness phenomenon
of nonlocal minimal surfaces have been recently performed in~\cite{MR3982031}.
\medskip

We also remark that nonlocal minimal surfaces
with a graphical structure enjoy special regularity properties:
for instance they are smooth inside their domain of minimization,
as established in~\cite{MR3934589}.

\section{Growth and Bern\v{s}te\u{\i}n properties
of minimal graphs}

A classical question, dating back to~\cite{zbMATH02583554} is whether or not minimizers of a given geometric
problem which have a global graphical structure are necessarily
hyperplanes. If this feature holds,
we say that the problem enjoys the ``Bern\v{s}te\u{\i}n property''.

This question is also related with the growth of these minimal
graphs at infinity. We now discuss and compare the classical
and the nonlocal worlds with respect to these features.

\subsection{Growth and Bern\v{s}te\u{\i}n properties
of classical minimal graphs}

The Bern\v{s}te\u{\i}n problem for classical minimal surfaces
asks whether or not solutions of
$$
{\rm div}\left(\frac{\nabla u}{\sqrt{1+|\nabla u|^2}}\right)=0 \;\; {\mbox{ in }}\R^n$$
are necessarily affine functions.

When~$n=2$, a positive answer to this question
was provided in~\cite{zbMATH02583554}.
The development of the theory related to the classical
Bern\v{s}te\u{\i}n problem is closely linked
to the regularity theory of classical minimal surfaces,
since, as proved in~\cite{MR178385},
the falsity of the Bern\v{s}te\u{\i}n property
in~$\R^n$ would imply the existence of a singular
minimal surface in~$\R^{n}$.
\medskip

In light of this connection,
the Bern\v{s}te\u{\i}n property for classical minimal surfaces
was established in~\cite{MR178385} for~$n=3$,
in~\cite{MR200816} for~$n=4$,
and in~\cite{MR233295} up to~$n=7$
(compare with the references in Section~\ref{INTCLA}).
A counterexample when~$n=8$ was then provided in~\cite{MR250205},
by constructing an eight-dimensional graph in~$\R^9$
which has vanishing mean curvature without being
affine.\medskip

As a counterpart of this kind of problems,
a precise growth analysis of classical minimal surfaces
with graphical structure has been obtained in~\cite{MR157096}
when~$n=2$ and in~\cite{MR248647} for all dimensions
(see also~\cite{MR296832, MR308945, MR412605, MR843597, MR2796515}
for related results).

More precisely, if~$u$ is a solution of
$$
{\rm div}\left(\frac{\nabla u}{\sqrt{1+|\nabla u|^2}}\right)=0 \;\; {\mbox{ in }}B_R\subset
\R^n,$$
then
\begin{equation} \label{GRO1}\sup_{B_R}|\nabla u|\le \exp\left[C\,\left(1+
\frac{\displaystyle\sup_{B_{2R} }u-u(0)}{R}\right)\right],\end{equation}
for some~$C>0$ depending only on~$n$.
Interestingly, this exponential type of gradient estimate is optimal,
as demonstrated in~\cite{MR157096}.

\subsection{Growth and Bern\v{s}te\u{\i}n properties
of nonlocal minimal graphs}

The Bern\v{s}te\u{\i}n properties
of nonlocal minimal graphs have been investigated
in~\cite{MR3680376, 2017arXiv170605701F}. In particular,
it is proved in these works that
the falsity of the nonlocal Bern\v{s}te\u{\i}n property
in~$\R^n$ would imply the existence of a singular nonlocal
minimal surface in~$\R^{n}$, thus providing a
perfect counterpart of~\cite{MR178385} to the nonlocal setting.\medskip

Combining this with the regularity results for nonlocal minimal surfaces
in~\cite{MR3090533}
and~\cite{MR3107529}, one obtains that the Bern\v{s}te\u{\i}n property
holds true for solutions~$u:\R^n\to\R$
of the $\sigma$-mean curvature equation in~\eqref{HS} provided that either~$n\in\{1,2\}$,
or~$n\le7$ and~$\sigma$ is sufficiently close to~$1$.
\medskip

It is an open problem to prove or disprove the
Bern\v{s}te\u{\i}n property for nonlocal minimal graphs~$u:\R^n\to\R$
in the general case
when~$n>2$ and~$\sigma\in(0,1)$
(the analysis on symmetric cones performed in~\cite{MR3798717} might
suggest that new phenomena could arise in dimension~$7$).\medskip

Concerning the growth of nonlocal minimal surfaces,
it has been established in~\cite{MR3934589}
that if~$u:B_{2R}\subset\R^n\to\R$ is a solution of the
$\sigma$-mean curvature equation in~$B_{2R}$,
then
\begin{equation} \label{GRO2} \sup_{B_R}|\nabla u|\le C\,\left(1+\frac{\displaystyle\sup_{B_R}u-\displaystyle\inf_{B_R}u}R\right)^{n+1+\sigma},\end{equation}
for a suitable constant~$C>0$,
depending only on~$n$ and~$\sigma$.

On the one hand, comparing~\eqref{GRO1} and~\eqref{GRO2}, we observe that
the nonlocal world offers us a wealth of surprises:
first of all, the estimate in~\eqref{GRO2}
is of polynomial type, rather than of exponential type, as it was
the bound in~\eqref{GRO1}. Roughly speaking,
this is due to the fact that two different and sufficiently
close``ends'' of a nonlocal minimal surface
repel each other (compare the discussion about the nonlocal
catenoid on page~\pageref{CATE}). Moreover, the right hand side
of~\eqref{GRO2} presents the oscillation of the solution in the
same ball~$B_R$, while~\eqref{GRO1} needed to extend the
right hand side to a larger ball such as~$B_{2R}$.

On the other hand, the nonlocal world maintains its own special difficulties:
for instance, differently from the classical case,
it is an open problem to determine whether the estimate in~\eqref{GRO2}
is optimal.\medskip

See Table~\ref{tab:1} for a sketchy description
of some similarities and differences between classical
and nonlocal minimal surfaces.

\begin{table*}[t]
\begin{center}{\large
  \begin{tabular}{ | p{3.5cm} | p{4.9cm} | p{4.9cm} | }
    \hline
     & {\em classical minimal surfaces} & {\em nonlocal minimal surfaces} \\ \hline
{\em interior regularity} & up to $n=7$ (optimal) & 
known when $n=2$, and up to $n=7$ if $\sigma$
is large enough (optimality unknown, no singular example available) \\ \hline
    {\em boundary regularity} & yes for convex (or mean convex) domains & 
boundary discontinuity (but enhanced differentiable regularity if continuous) \\    \hline
{\em Bern\v{s}te\u{\i}n property for~$u:\R^n\to\R$}
& up to $n=7$ (optimal) & known when $n=2$, and up to $n=7$ if $\sigma$
is large enough (optimality unknown, no counterexample available) \\    \hline
{\em gradient growth at infinity} & exponential (optimal) & polynomial (optimality unknown) \\    \hline
  \end{tabular}}
  \caption{Classical versus nonlocal minimal surfaces}
  \label{tab:1}
\end{center}\end{table*}

\section{Phase coexistence models}

A physical situation in which minimal interfaces naturally
arise occurs in the mathematical description of phase
coexistence (or phase transition). In the classical setting,
the ansatz is that the short range particle interaction produces
a surface tension that leads, on a large scale,  to interfaces
of minimal area. A modern variation of this model takes into
account long-range particle interactions, and we describe
here in which sense this model is related to the theory of (non)local
minimal surfaces.

\subsection{Classical phase coexistence models}\label{12345asasCL}

One of the most popular phase coexistence models
is the one producing the so-called Allen-Cahn equation
\begin{equation}\label{ALLE} -\Delta u=u-u^3.\end{equation}
Equation~\eqref{ALLE} has a variational structure, corresponding
to critical points of an energy functional of the type
\begin{equation}\label{ALLE1} 
{\mathcal{F}}(u;\Omega)=
{\mathcal{F}}(u):=\int_\Omega |\nabla u(x)|^2+
W(u(x))\,dx,\end{equation}
where~$W$ is a ``double-well'' potential attaining its minima
at the ``pure phases'' $-1$ and~$+1$.

A classical link between the functional in~\eqref{ALLE1} and
the perimeter can be expressed in terms of~$\Gamma$-convergence,
as stated in~\cite{MR445362}. More precisely,
given a smooth domain~$\Omega\subset\R^n$ and~$\e>0$, one can
consider the rescaled functional
\begin{equation}\label{ugq34} {\mathcal{F}}_\e(u;\Omega)={\mathcal{F}}_\e(u):=\int_\Omega \e\,|\nabla u(x)|^2+\frac1\e\,
W(u(x))\,dx.\end{equation}
Interestingly, the functional~${\mathcal{F}}_\e$
in~\eqref{ugq34} inherits the minimization properties of
the functional~${\mathcal{F}}$
in~\eqref{ALLE1} after a dilation: for instance,
if~$u$ is a minimizer for~${\mathcal{F}}$
in every ball, then
\begin{equation}\label{uep}u_\e(x) := u\left(\frac{x}\e\right)\end{equation}
is a minimizer for~${\mathcal{F}}_\e$
in every ball.
Moreover, it holds that:
\begin{itemize}
\item for any $u_\e:\R^n\to[-1,1]$ converging to~$u$
in~$L^1_{\rm loc}(\R^n)$, we have that
\begin{equation}\label{IS-1239}\liminf_{\e\searrow0}{\mathcal{F}}_\e(u_\e)\ge\mathcal{F}_0(u),\end{equation}
\item given~$u:\R^n\to[-1,1]$, there exists~$u_\e:\R^n\to[-1,1]$converging to~$u$ in~$L^1_{\rm loc}(\R^n)$, such that
\begin{equation}\label{IS-1239-2}\limsup_{\e\searrow0}{\mathcal{F}}_\e(u_\e)\le\mathcal{F}_0(u),\end{equation}
\end{itemize}
where, for a suitable~$c>0$,
\begin{equation}\label{F0}{\mathcal{F}}_0(u;\Omega)=
\mathcal{F}_0(u):=
\begin{cases}
c\,{\rm Per}(E,\Omega) &
\begin{matrix}
{\mbox{if $u=\Xi_E$}} \\{\mbox{for some set~$E$,}}
\end{matrix}\\
\\
+\infty &{\mbox{otherwise.}}\end{cases}\end{equation}
This type of ``functional convergence''
has also a ``geometric counterpart'', given by the locally
uniform convergence of levels sets, as established in~\cite{MR1310848}.
More precisely, if~$u$ is a minimizer for~${\mathcal{F}}$
in every ball and~$u_\e$ is as in~\eqref{uep},
then~$u_\e$ converges, up to subsequences, in~$L^1_{\rm loc}(\R^n)$
to some~$u= \Xi_E$, with~$E$ minimizing the perimeter,
and the level sets of~$u_\e$ approach~$\partial E$ locally uniformly,
that is, given any~$\vartheta\in(0,1)$,
any~$R > 0$ and any~$\delta > 0$
there exists~$\e_0\in(0,1)$
such that, if~$\e\in(0,\e_0)$, we have that
\begin{equation}\label{LSC}
\big\{ |u_\e|<1-\theta\big\}\cap B_R\subseteq\bigcup_{p\in\partial E}
B_\delta(p).
\end{equation}
The level set convergence in~\eqref{LSC}
relies on suitable
``energy and density estimates'' stating that 
the energy behaves as an ``interface'' and,
in a neighborhood of
the interface, both the phases have positive densities. More precisely,
if~$R\ge1$ and~$u$ is a minimizer of~${\mathcal{F}}$
in~$B_{R+1}$, then
\begin{equation}\label{ENDH:es} {\mathcal{F}}_{B_R}(u)\le CR^{n-1},\end{equation}
for some~$C>0$, and, for every~$\vartheta_1,\vartheta_2\in(0,1)$,
if~$|u(0)|\le\vartheta_1$ then
\begin{equation}\label{DENISUI}
\begin{split}&
\{u\ge\vartheta_2\}\cap B_R\ge cR^n \\{\mbox{and}}\qquad&
\{u\le-\vartheta_2\}\cap B_R\ge cR^n.\end{split}
\end{equation}

\subsection{Long-range phase coexistence models}

In the recent literature a number of models
have been introduced in order to study phase coexistence
driven by long-range particle interactions, see e.g.~\cite{MR1612250, MR2223366}.

We describe here a simple model
of this kind, comparing the results
obtained in this situation with the classical ones presented in Section~\ref{12345asasCL}.\medskip

We focus on a nonlocal version of the Allen-Cahn equation~\eqref{ALLE}
given by
\begin{equation}\label{ALLEfrac} (-\Delta)^s u=u-u^3,\end{equation}
with~$s\in(0,1)$.
Equation~\eqref{ALLEfrac} has a variational structure, corresponding
to critical points of an energy functional of the type
\begin{equation}\label{ALLE1frac} \begin{split}&
{\mathcal{E}}(u;\Omega)=
{\mathcal{E}}(u):=\int_{Q_\Omega} \frac{|u(x)-u(y)|^2}{|x-y|^{n+2s}}\,
dx\,dy\\&\qquad\qquad+\int_\Omega
W(u(x))\,dx,\end{split}\end{equation}
where~$W$ is a ``double-well'' potential attaining its minima
at the ``pure phases'' $-1$ and~$+1$, and~$Q_\Omega$
is the cross-type domain given by
$$ (\Omega\times\Omega)\cup(\Omega\times\Omega^c)\cup
(\Omega^c\times\Omega).$$
The rationale behind this choice of~$Q_\Omega$
is to collect all the couples of points in which at least
one of the two lies in~$\Omega$, so to comprise all the
point interactions that involve~$\Omega$ (a similar choice
was performed in the geometric problem~\eqref{PERs3}).

The $\Gamma$-convergence theory for the functional in~\eqref{ALLE1frac}
has been established in~\cite{MR2948285}.
For this, one needs to introduce a family of rescaled functional
with different scaling properties depending on the fractional parameter~$s$.
More precisely, it is convenient to define
\begin{equation} \label{Eeps}\begin{split}&
{\mathcal{E}}_\e(u;\Omega)=
{\mathcal{E}}_\e(u):=a_\e\,\int_{Q_\Omega} \frac{|u(x)-u(y)|^2}{|x-y|^{n+2s}}\,
dx\,dy\\&\qquad\qquad+b_\e\,\int_\Omega
W(u(x))\,dx,\end{split}\end{equation}
where
$$ a_\e:=\begin{cases}
\e^{2s-1} & {\mbox{ if }}s\in(1/2,\,1),\\
|\log\e|^{-1} & {\mbox{ if }}s=1/2,\\
1& {\mbox{ if }}s\in(0,\,1/2),
\end{cases}$$
and
$$ b_\e:=\e^{-2s}a_\e=\begin{cases}
\e^{-1} & {\mbox{ if }}s\in(1/2,\,1),\\
\e^{-1}|\log\e|^{-1} & {\mbox{ if }}s=1/2,\\
\e^{-2s}& {\mbox{ if }}s\in(0,\,1/2).
\end{cases}$$
The coefficients~$a_\e$ and~$b_\e$ are tuned with
respect to the energy of a one-dimensional layer joining~$-1$
and~$+1$ and, in this way, the functional~${\mathcal{E}}_\eps$
in~\eqref{Eeps} satisfies some interesting $\Gamma$-convergence properties.
Remarkably, these properties are very sensitive to the fractional
parameter~$s$: indeed, 
when~$s<1/2$ the $\Gamma$-convergence theory
obtains the $\sigma$-perimeter introduced in~\eqref{PERs3},
with~$\sigma:=2s$,
but when~$s\ge1/2$, the nonlocal
character of the problem is lost in the limit and the
$\Gamma$-convergence theory boils down to the classical one in~\eqref{IS-1239}, \eqref{IS-1239-2} and~\eqref{F0}.

More precisely, one can define
\begin{equation*}{\mathcal{E}}_0(u;\Omega)=
\mathcal{E}_0(u):=
\begin{cases}
c\,{\rm Per}(E,\Omega) &
\begin{matrix}
{\mbox{if $s\in[1/2,1)$}}\\{\mbox{and $u=\Xi_E$}} \\{\mbox{for some set~$E$,}}
\end{matrix}\\
\\
c\,{\rm Per}_{2s}(E,\Omega) &
\begin{matrix}
{\mbox{if $s\in(0,1/2)$}} \\{\mbox{and $u=\Xi_E$}} \\{\mbox{for some set~$E$,}}
\end{matrix}
\\
\\
+\infty &{\mbox{otherwise.}}\end{cases}\end{equation*}
In this sense, the functional~${\mathcal{E}}_0$ is the ``natural
replacement'' of the one in~\eqref{F0} to deal with the~$\Gamma$-convergence
of nonlocal phase transitions, in the sense that the statements
in~\eqref{IS-1239}, \eqref{IS-1239-2} and~\eqref{F0} hold true in this case with~${\mathcal{F}}_\e$
replaced by~${\mathcal{E}}_\e$ and~${\mathcal{F}}_0$
replaced by~${\mathcal{E}}_0$.\medskip

For additional discussions of nonlocal $\Gamma$-convergence results,
see~\cite{2020arXiv200101475D}.\medskip

The convergence of level sets in~\eqref{LSC}
has also a perfect counterpart for the minimizers of~${\mathcal{E}}$,
as established in~\cite{MR3133422}.
These results also rely on suitable density estimates
which guarantee that~\eqref{DENISUI} also holds true for
minimizers of~${\mathcal{E}}$.\medskip

However, in this case, an energy estimate as in~\eqref{ENDH:es}
exhibits a significant difference, depending on the fractional
exponent~$s$. Indeed, in this framework, we have that
if~$R\ge2$ and~$u$ is a minimizer of~${\mathcal{E}}$
in~$B_{R+1}$, then
\begin{equation} \label{Olsdepr}{\mathcal{E}}_{B_R}(u)\le \begin{cases}
CR^{n-1} & {\mbox{ if }}s\in(1/2,\,1),\\
CR\,\log R& {\mbox{ if }}s=1/2,\\
CR^{n-2s} & {\mbox{ if }}s\in(0,\,1/2),
\end{cases}\end{equation}
for some~$C>0$.

In spite of its quantitative difference with~\eqref{ENDH:es},
the estimate in~\eqref{Olsdepr} is sufficient to determine that
the interface contribution is ``negligible'' with respect to the Lebesgue
measure of large balls. Moreover, such an estimate is in agreement
with the energy of one-dimensional layers.\medskip

We refer to Table~\ref{tab:2} for a schematic representation
of similarities and differences of classical and nonlocal
phase transitions in terms of $\Gamma$-limits and density estimates.\medskip

In addition, long- phase transition models and nonlocal
minimal surfaces have close connections with
spin models and discrete systems, see~\cite{MR1372427, MR1453735, MR1794850, MR2219278, MR2608835, MR3017265, MR3472011} and the references therein.
In particular, in~\cite{MR3652519} the
reciprocal approximation of ground states for long-range Ising models and nonlocal minimal surfaces
has been established and the relation between these two
problems has been investigated in detail.

\begin{table*}[t]
\begin{center}{\large
  \begin{tabular}{ | l | p {4 cm} | p {4 cm} | p {4 cm} |}
    \hline
     & {\em phase transitions $s\in(0,1/2)$} & {\em phase transitions $s=1/2$} & {\em phase transitions $s\in(1/2,1]$}\\ \hline
    {\em $\Gamma$-limit} & nonlocal perimeter & classical perimeter & classical perimeter \\ \hline
    {\em energy estimates} & $R^{n-2s}$ &$ R\,\log R $&$ R^{n-1}$ \\    \hline
  \end{tabular}}
  \caption{Classical versus nonlocal phase transitions}
  \label{tab:2}
\end{center}\end{table*}

\section{Geometric flows}

Another classical topic in geometric analysis consists in the
study of evolution equations describing the motion of set boundaries.
In this framework, at any time~$t$, one prescribes
the speed~$V$ of a hypersurface~$\partial E_t$ in the direction of the outer normal to~$E_t$.

In this section, in light of the setting in~\eqref{HS},
we will describe analogies and differences between
the classical and the nonlocal mean curvature flows.

\subsection{Mean curvature flow}\label{HANSX:SSKS}

One of the widest problem in geometric analysis
focuses on the case in which
the normal velocity~$V$ of the flow is taken to
be the classical mean curvature. This flow
decreases the classical perimeter and has a number of remarkable
properties.
A convenient description of the mean curvature flow can be given
in terms of viscosity solutions, using a level set approach,
see~\cite{MR1100206}. See also~\cite{MR1100211, MR1360824} for related
viscosity methods. 
If the initial hypersurface satisfies a
Lipschitz condition, then the mean curvature flow possesses a short time existence theory, see~\cite{MR1117150}.
In general, the flow may develop
singularities, and the understanding of the different possible
singularities is a very rich and interesting field of investigation in itself,
see e.g.~\cite{MR485012, MR2024995, MR2815949}.
\medskip

The planar case~$n=2$ presents some special feature. In this case,
the flow becomes ``curve shortening'' in the sense that
its special characteristic is to decrease the length of a given closed curve
as fast as possible.
In this setting, it was proved in~\cite{MR0906392} that the
curve shortening flow makes
every smooth closed embedded curves in the plane shrink smoothly 
to round points (i.e., the curve shrinks becoming ``closer and closer to
a circle''). More precisely: on the one hand, as shown in~\cite{MR0840401},
if the initial curve is convex, then the evolving curve
remains convex and becomes asymptotically
circular as it shrinks; on the other hand, as shown in~\cite{MR0906392},
if the curve is only embedded, then
no singularity develops before the curve becomes convex, then shrinking to
a round point. 

These results can be also interpreted in light of
a ``distance comparison property''. Namely, given two points~$p$
and~$q$ on the evolving surface, one can consider
the ``extrinsic'' distance~$d=d(p,q)$, given by the Euclidean distance
of~$p$ and~$q$, and the arc-length ``intrinsic'' distance~$\ell=\ell(p,q)$,
that measures the distance of~$p$ and~$q$ along the curve.
For a closed curve of length~$L$, one can also define
$$ \lambda:=\frac{L}\pi\,\sin\frac{\pi\ell}{L}.$$
Then, as proved in~\cite{MR1656553}, we have that the minimal ratio~$\rho:=\min\frac{d}\lambda$  is nondecreasing under the curve shortening flow
(and, in fact, strictly increasing unless the curve is a circle).
In this spirit, the ratio~$\rho$ plays the role of
an improving isoperimetric quantity which
measures the deviation of the evolving curve from a round circle.
In particular, if a curve developed a ``neckpinch'', it would present two
points for which~$d=0$ and~$\ell>0$, hence~$\rho=0$: as a consequence,
the monotonicity proved in~\cite{MR1656553} immediately implies
that planar curves do not develop 
neckpinches under the curve shortening flow. More generally,
the results in~\cite{MR1656553} can be seen as a quantitative
refinement of those in~\cite{MR0840401, MR0906392}.
See also~\cite{MR0772132} for related results.

In any case, for the mean curvature flow, the formation of
neckpinch singularity for the evolving hypersurface can only occur
(and it occurs) in dimension~$n\ge3$, see~\cite{MR1030675, MR3803553}.

\subsection{Nonlocal mean curvature flow}

A natural variant of the flow discussed in Section~\ref{HANSX:SSKS}
consists in taking as normal velocity the nonlocal mean curvature of
the evolving set, as given in~\eqref{HS}. The main feature of
this evolution problem is that it
decreases the nonlocal perimeter defined in~\eqref{PERs3}.

This evolution problem emerged in~\cite{MR2487027, MR2564467}.
A detailed viscosity solution theory for the nonlocal mean curvature
flow has been presented in~\cite{MR3156889}, see also~\cite{MR3401008}
for more general results.

In the viscosity setting, given an initial set~$E\subset\R^n$, 
the evolution of~$E$ is described by level sets of a function~$u_E=u_E(x,t)$
which solves a parabolic problem driven by
the nonlocal  mean curvature. More precisely, the evolution of~$E$
is trapped between an outer and an inner flows defined by
\begin{eqnarray*}&& E^+(t):=\{x\in\R^n {\mbox{ s.t. }}u_E(x,t)\ge0\}
\\{\mbox{and }}&&
E^-(t):=\{x\in\R^n {\mbox{ s.t. }}u_E(x,t)>0\}.\end{eqnarray*}
The ideal case would be the one in which these two flows coincide
``up to a smooth manifold'' which describes the evolving surface,
i.e. the case in which
\begin{equation}\label{HAJS:38i4r345y77} \begin{split}
\Sigma_E(t)\,&:=\,E^+(t)\setminus E^-(t)\\&
=\,
\{x\in\R^n {\mbox{ s.t. }}u_E(x,t)=0\}\end{split}\end{equation}
is a nice hypersurface, but, in general, such a property is not
warranted by the notion of viscosity solutions.\medskip

A short-time existence theory of smooth solutions of the nonlocal mean curvature flow has been recently 
established in~\cite{2019arXiv190610990J}
under the assumption that the initial hypersurface is of
class~$C^{1,1}$. As far as we know,
it is still an open problem to determine whether
the nonlocal mean curvature flow possesses 
a short-time existence theory of smooth solutions for Lipschitz initial
data (this would provide a complete
counterpart of the classical results in~\cite{MR1117150}). 
\medskip

It is certainly interesting to detect suitable information
that are preserved by the nonlocal mean curvature flow.
For instance:
the fact that the nonlocal mean curvature flow preserves the positivity
of the nonlocal mean curvature
itself has been established in~\cite{MR3951024},
and the preservation of convexity has been proved in~\cite{MR3713894}.
\medskip

An important difference between the classical and the nonlocal
mean curvature flows consists in the ``fattening phenomena''
of viscosity solutions, in which the set in~\eqref{HAJS:38i4r345y77}
may develop a nonempty interior (thus failing to be a nice
hypersurface),
as investigated in~\cite{MR4000255}. For instance,
one can consider the case in which the initial set is the ``cross in the plane''
given by
\begin{equation}\label{E7js023d4} E=\{ (x,y)\in\R^2{\mbox{ s.t. }}|x|>|y|\}.\end{equation}
In this situation, 
using the notation in~\eqref{HAJS:38i4r345y77},
one can prove that~$\Sigma_E(t)$ is nontrivial
and immediately develops a positive measure: more precisely,
we have that~$\Sigma_E(t)$ contains the ball~$ B_{ct^{\frac1{1+s}}}$,
for some~$c>0$.\medskip

On the one hand,
this result is in agreement with the classical
case, since also the classical mean curvature flow develops
this kind of situations when starting from the set in~\eqref{E7js023d4}
(see~\cite{MR1100206}).\medskip

On the other hand, the fattening phenomena
for the nonlocal situation offers a number
of surprises with respect to the classical case (see~\cite{MR4000255}).
First of all, the quantitative properties
of the interaction kernels play a decisive role
in the development of the fat portions of~$\Sigma_E(t)$
and in general on the evolution of the set~$E$. For instance,
if the interaction kernel~$\frac{1}{|x-y|^{n+\sigma}}$ in~\eqref{IN3r45T} is replaced by
a different kernel that is smooth and compactly supported,
then the nonlocal mean curvature flow of the cross~$E$ in~\eqref{E7js023d4}
does not develop any fattening and, more precisely,
the evolution of~$E$ would be~$E$ itself (in a sense,
we switch from a fattening situation, produced
by singular kernels with slow decay, to a ``pinning effect''
produced by smooth and compactly supported kernels).\medskip

Furthermore, there are cases in which the fattening
phenomena are different according to the type of local or nonlocal
mean curvature flow that we consider (see~\cite{MR1100206}). For instance,
if the initial set is made by two tangent balls in the plane such as
\begin{equation}
E:= B_1(-1,0)\cup B_1(1,0),
\end{equation}
then, for the nonlocal mean curvature flow,
using the notation in~\eqref{HAJS:38i4r345y77}, we have that~$\Sigma_E(t)$
has empty interior for all~$t>0$.
This situation is different with respect to the classical
mean curvature flow, which immediately develops
fattening (see~\cite{MR1298266}).
\medskip

The nonlocal mean curvature flow (and its many variants)
are a great source of intriguing questions and mathematical
adventures. Among the interesting questions remained open,
we mention that
it is not known whether the nonlocal mean curvature flow
possesses a nonlocal version of the monotonicity formula in~\cite{MR1656553},
or what a natural replacement of it could be. This question
is also related to a deeper understanding of the singularities
of the nonlocal mean curvature flow.

\section*{Acknowledgments}

The author is member of INdAM/GNAMPA and AustMS.
Her research is supported by the Australian Research Council Discovery Project DP170104880 ``N.E.W. Nonlocal Equations at Work''
and by the DECRA Project DE180100957 ``PDEs, free boundaries and applications''.

It is a pleasure to thank Enrico Valdinoci for his comments
on a preliminary version of this manuscript.

\end{document}